# From The Goldbach Conjecture To The Theorem

**Pereyra, P.H.**[1]
UNISC-FACCAT-FAPA

pereyraph@gmail.com
www.pereyra.hostmach.com.br

**Bodmann, B.E.J**
UFRGS

bejbodmann@gmail.com

### Abstract

In the present work we demonstrate that the so called Goldbach conjecture from 1742 -- "All positive even numbers greater than two can be expressed as a sum of two primes" – due to Leonhard Euler, is a true statement. This result is partially based on the Wilson theorem, and complementary on our reasoning to cast the problem into a diophantine equation. The latter is the master equation for the conjectures proof.

*(1) Corresponding Author*

## INTRODUCTION

The Goldbach conjecture is probably one of the most enduring challenges in mathematics, asserting that every even number larger than two may be expressed as the sum of two primes. Some progress has been made on related weaker assertions, there also exist numerical confirmations of the conjecture up to huge numbers. The conjecture as it stands shall hold for all even natural numbers larger than 4, which means that numerical proofs are valid only for a finite dimensional subset. For some attempts on this subject see for instance Schroeder (1997) [1] and references therein.

## THEOREM

All positive even integer numbers greater than two can be expressed as a sum of two primes.

## PROOF

As a first step we consider the trivial fact that any even number $2n$ $(n \in N, n > 2)$ may be expressed as the sum of two odd numbers. An illustrative example for the number $32$ is given in figure 1, where we introduce a distribution matrix. Here the row and column represent the odd summands resulting in $2n$, while the matrix elements show the possible sums.

Fig.1 - DISTRIBUTION MATRIX

One may conclude by induction that an even number $2n$ can be represented as a sum of two odd numbers from exactly $|I_n|$ combinations given by the size

$$|I_n| = \frac{n - 2 + (n \bmod 2)}{2}, \qquad (1)$$

where each number is the sum of the $i\_th$ element of the horizontal set of odd numbers, defining the columns of the distribution matrix, and given by sequence

$$I_n = \{2n-3, \ldots, 2n-1-2|I_n|\}, \qquad (2)$$

with the i_th element of the vertical set (interval) of odd numbers, defining respectively the rows of the distribution matrix, given by sequence

$$\bar{I}_n = \{3, \ldots, 1+2|I_n|\}. \qquad (3)$$

The sum is then

$$I_n^i + \bar{I}_n^i = 2n, \forall n, \forall i, \{i \in N \mid 1 \leq i \leq |I_n|\}. \qquad (4)$$

In the afore mentioned example the number *32* is represented by the sum of two odd numbers in exactly *7* different ways.

A further restriction is introduced by transforming the distribution matrix into a boolean matrix, shown in figure 2. The horizontal and vertical prime numbers are indicated by 1, whereas all remaining composed numbers are tagged by 0. The fields of the distribution matrix contains even numbers, which are genuine compositions of two primes and are indicated by 1, whereas the remainders, i.e. sum of at least one non-prime, are represented by 0.

Fig.2- BOOLEAN MATRIX

By observation one verifies that in the diagonal line of the number *32* there appear two true combinations (marked with a box), indicating that the number *32* can be formed by the sum of two prime numbers only in two ways (*29+3=32* and *19+13=32*).

Considering the Wilson Theorem [2] and the consequent binary formula [4],[6] for primes given by the floor function $\lfloor x \rfloor$

$$F(j) = \left\lfloor \cos^2\left(\pi \frac{(j-1)!+1}{j}\right) \right\rfloor, \qquad (5)$$

we can express the corresponding boolean sets (2) and (3) for the number *2n* as

$$B_n = [F(2n-3), \ldots, F(2n-1-2|I_n|)] \qquad (6)$$

and

$$\overline{B}_n = [F(3), \ldots, F(1+2|\overline{I}_n|)]. \qquad (7)$$

Consequently the number of representations as sums of two prime numbers for the even number *2n* is given by the dot product of vector like forms (6) and (7)

$$|B_n| = B_n \bullet \overline{B}_n = \sum_{i=1}^{|I_n|} B_n^i \overline{B}_n^i, \qquad (8)$$

where $B_n^i$ and $\overline{B}_n^i$ correspond to the *i_th* elements of the respective sets (6) and (7).

For the proof of the theorem it is sufficient to demonstrate that

$$|B_n| \geq 1, \quad \forall n, \; n \in N \mid n > 2 \qquad (9)$$

and because of (5) and (8) one concludes that (9) is valid when there exists at least one term in (8) with

$$\lfloor \cos^2(\alpha(I_n^i)\pi) \rfloor \lfloor \cos^2(\alpha(\overline{I}_n^i)\pi) \rfloor = 1 \Leftrightarrow \cos(\alpha(I_n^i)\pi)\cos(\alpha(\overline{I}_n^i)\pi) = \pm 1, \qquad (10)$$

where $I_n^i$ e $\overline{I}_n^i$ correspond to the *i-th* elements of the sequences (2) and (3), respectively, and

$$\alpha(I_n^i) = \frac{(I_n^i - 1)! + 1}{I_n^i} \qquad (11)$$

is called phase function.

From the Wilson theorem [2] follows that both, $I_n^i$ and $\bar{I}_n^i$ are prime numbers if $\alpha(I_n^i) = r$ and $\alpha(\bar{I}_n^i) = s$ with $r, s \in N^*$, so that condition (10) shall occur for some $i \in N \mid 1 \leq i \leq |I_n|$.

From (2) and (3) we see that the elements $I_n^i$ e $\bar{I}_n^i$ possess the following relation as a function of $i$

$$\begin{aligned}\bar{I}_n^i &= 2i+1 \\ I_n^i &= 2(n-1-i)+1\end{aligned} \quad i = 1..|I_n| \quad , \tag{12}$$

then if $j=2i+1$, (8) represents the discrete convolution of (5) in the point $n-1$ given by

$$|B_n| = B_n \bullet \bar{B}_n = \sum_{i=1}^{|I_n|} B_n^i \bar{B}_n^i = [F(j) \otimes F(j)](n-1) = \sum_{i=1}^{|I_n|} F(2i+1)F(2(n-1-i)+1) \ . \tag{13}$$

If (13) implies (9) this proofs the theorem.

It is instructive to use an analogy by wave superposition and observe that equation (13) is equivalent to two irregular waves having some sort of maximum constructive interference,

$$\cos^2\left(\pi \frac{(2i)!+1}{2i+1}\right) \cos^2\left(\pi \frac{(2(n-1-i))!+1}{2(n-1-i)+1}\right) = 1 \ , \text{ for some natural } i \ (i \in N \mid 1 \leq i \leq |I_n|) . \tag{14}$$

In other words, it is necessary to demonstrate by (14), that there always exists at least one maximum for some natural $i$ and natural $n$.

From the Wilson theorem [2] we know that for $I_n^i, \bar{I}_n^i$ primes follows

$$\begin{aligned}(2i)!+1 &= 0 \mod (2i+1) \Rightarrow \frac{[(2i)!+1]}{(2i+1)} = a \\ (2n-2-2i)!+1 &= 0 \mod (2n-1-2i) \Rightarrow \frac{[(2n-2-2i)!+1]}{(2n-1-2i)} = b\end{aligned} \quad , a,b \in N^* , \tag{15}$$

and equaling the terms $2i+1$ in (15), implies in a non linear diophantine equation

$$b[(2i)!+1] + a[(2n-2-2i)!+1] = 2nab \ . \tag{16}$$

## ZERO APLICATION - SYNCHRONIZATION PROOF

We can write (16) as

$$b[(2i)!+1] + a[(2i)!+1+(2n-2-2i)!+1-(2i)!-1] = 2nab, \qquad (17)$$

and considering $a = \alpha$ a natural number and $2i+1$ a prime number for some $i \in N \mid 1 \le i \le |I_n|$, such that $\alpha(2i+1) = (2i)!+1$, as guaranteed from the Wilson theorem [2], Chebyshev [4] theorem, and from the inequality bounds by the prime counting function of Rosser and Schoenfeld [3], yields

$$b(2n-1-2i) - \alpha(2i+1) = (2n-2-2i)! - (2i)! \, , \, \alpha \in N^*, \qquad (18)$$

that is a linear diophantine equation with determinate $\alpha$ and $i$.

An equivalent form of (18) is

$$b(2i+1)\left(\frac{2n}{(2i+1)} - 1\right) - \alpha(2i+1) = (2n-2-2i)! - (2i)! \, , \, \alpha \in N^* \qquad (19)$$

and follows that if some $2i+1$ prime number in the interval $i \in N \mid 1 \le i \le |I_n|$ does not divide $n$, implies that $\gcd\left((2i+1)\left(\frac{2n}{(2i+1)} - 1\right), -(2i+1)\right) = 1$, and from the theorem of solutions of the linear diophantine equations [5], (19) and consequently (18) have natural solutions,. This necessarily happens for some prime number such that $2i+1 > \frac{n}{2}$, because of the inequality bounds by the prime counting function of Rosser and Schoenfeld [3] and is verified directly for small $2n$. This not exclude the possibility of the existence of the natural solutions for (18) with $n$ multiple of the $2i+1$.

From (19), we can change the coefficient of the $(2i+1)$ by $(-\alpha+k)$ units, where $k \in N^*$ is an arbitrary natural number, then this results

$$b(2i+1)\left(\frac{2n}{(2i+1)} - 1\right) - (\alpha - \alpha + k)(2i+1) = (2n-2-2i)!+1 - k(2i+1) \, , \, \alpha, k \in N^* \, , \qquad (20)$$

and simplifying yields

$$b(2i+1)\left(\frac{2n}{(2i+1)} - 1\right) - k(2i+1) = (2n-2-2i)!+1 - k(2i+1) \, , \, k \in N^*. \qquad (21)$$

But, as $\gcd\left((2i+1)\left(\dfrac{2n}{(2i+1)}-1\right), -(2i+1)\right)=1$, this imply that (21) have natural solutions $(k, b_0)$, (note that $k$ not change because is an arbitrary natural number!), such that we can write

$$b_0(2i+1)\left(\dfrac{2n}{(2i+1)}-1\right)-k(2i+1) = (2n-2-2i)!+1-k(2i+1), \quad k, b_0 \in N^* \tag{22}$$

and (22) yields

$$b_0 = \dfrac{(2n-2-2i)!+1}{(2n-1-2i)}, \quad b_0 \in N^*, \tag{23}$$

then follow by Wilson theorem that $2n-1-2i$ also is a prime number.

This result show that the condition $\gcd\left((2i+1)\left(\dfrac{2n}{(2i+1)}-1\right), -(2i+1)\right)=1$ with $2i+1$ prime number for (18), always implies in the possibility of the synchronization with a $2n-1-2i$ prime number.

**THE CONVERSE PROOF**

In the further we show that (18) always posses solutions with $b=b_0$ naturals (non zero) numbers, and by inverse algebra the same imply in the master diophantine equation (16) and the Wilson conditions (15).

From (22), we can change the coefficient of the $(2i+1)$ by $\alpha$ units, and write

$$b_0(2i+1)\left(\dfrac{2n}{(2i+1)}-1\right)-\alpha(2i+1) = (2n-2-2i)!-(2i)!, \quad \alpha, b_0 \in N^* \tag{24}$$

showing that $(\alpha, b_0)$ is a natural solution of the (18), and as $\alpha$ is previously determinate in (17) also is solution for (16) by inverse algebra, with $2n-1-2i$ prime number. Then, consequently the Wilson conditions (15) are valid. This complete the converse proof.

We can say that (18) represent a difference between multiple of the two co-prime numbers, such that one of them is a prime number, in this case $2i+1$. Follow that if (18) have natural solutions, then obligatorily exist $2n-1-2i$ prime number, such that (18) also have natural solution.

## CONCLUSION

The present outline convinced us to have shown that Goldbach's conjecture is a genuine theorem. From the developments above one recognizes that there always exist a natural *a* and *b* that synchronize two primes $2i+1$ and $2n-1-2i$ in the interval $i \in N \mid 1 \leq i \leq |I_n|$ and its sum equal *2n*. Then (16) always has natural solutions in the interval $i \in N \mid 1 \leq i \leq |I_n|$ and consequently (14) is valid for any natural *n*. This implies that (9) is true, *quod erat demonstrandum*. The case *n=2* is verified directly.

Concluding, it is noteworthy that the development above in its general form does not exclude non integer solutions that imply in the representation of *2n* by a sum of two non prime numbers. In order to obtain this case it is sufficient to note that permits a rational *a* and/or *b*. Thus the present development shows that there exist three ways to represent the number *2n*. (A) As a sum of two prime numbers (with natural *a* and *b*), (B) as a sum of one prime number and an odd non prime number (with either *a* or *b* natural and the seconds rational), and (C) as a sum of two odd non prime numbers (with *a* and *b* rational). We see that the cases (B) and (C) depend on the existence of odd non prime numbers in the sequences (2) and (3). The case (A) is always feasible due to the existence of prime numbers in the sequences (2) and (3) as guaranteed by Chebyshev's theorem [4] and the inequality bounds by the prime counting function of Rosser and Schoenfeld [3].